\pdfoutput=1 
\documentclass[11pt,letterpaper]{amsart}
\PassOptionsToPackage{obeyspaces}{url}
\usepackage{amsfonts, amsmath, amssymb}

\usepackage[hidelinks,hypertexnames=false]{hyperref}
\usepackage{bookmark}
\bookmarksetup{open,numbered,depth=5} 

\usepackage{amsthm}

\usepackage{enumitem}

\usepackage{breakurl}
\usepackage{tikz}

\usepackage{etoolbox} 
\patchcmd{\thebibliography}{\leftmargin\labelwidth}{\leftmargin\labelwidth\addtolength\itemsep{-0.1\baselineskip}}{}{}

\usepackage{mathtools}

\usepackage{microtype}

\author{Boris Bukh}
\thanks{Supported in part by U.S.\ taxpayers through NSF CAREER grant DMS-1555149.}

\title{Extremal graphs without exponentially-small bicliques}

\usepackage[nameinlink]{cleveref}

\newtheorem{theorem}{Theorem}
\newtheorem{lemma}[theorem]{Lemma}

\newtheorem{proposition}[theorem]{Proposition}

\theoremstyle{definition}
\newtheorem{definition}[theorem]{Definition}

\AtBeginEnvironment{theorem}{\let\goodoldlabel\label \def\label{\goodoldlabel[theorem]}}
\AtBeginEnvironment{lemma}{\let\goodoldlabel\label \def\label{\goodoldlabel[lemma]}}
\AtBeginEnvironment{proposition}{\let\goodoldlabel\label \def\label{\goodoldlabel[proposition]}}

\newcommand*{\eqdef}{\stackrel{\mbox{\normalfont\tiny def}}{=}}  
\newcommand*{\veps}{\varepsilon}                                 
\DeclarePairedDelimiter\abs{\lvert}{\rvert}                      
\newcommand*{\Fq}{\mathbb{F}_q}                                  
\newcommand*{\Fqbar}{\overline{\mathbb{F}_q}}                    
\newcommand*{\Znonneg}{\mathbb{Z}_{\geq 0}}                      
\newcommand*{\E}{\mathbb{E}}                                     
\newcommand*{\V}{\mathbf{V}}                                     
\renewcommand*{\P}{\mathbb{P}}                                   
\DeclareMathOperator{\Var}{Var}                                  
\DeclareMathOperator{\Gr}{Gr}                                    
\DeclareMathOperator{\ex}{ex}                                    
\DeclareMathOperator{\vecspan}{span}                             
\newcommand*{\VD}{\Phi}
\newcommand*{\vd}{\phi}
\makeatletter
  \def\poly{\@ifnextchar[{\@polygrade}{\@polyplain}}
  \def\@polygrade[#1]#2#3#4{\@polyplain{#2}{#3}{#4}_{#1}}
  \def\@polyplain#1#2#3{#1[#2_0,\dotsc,#2_{#3}]}
\makeatother
  
\setlist[enumerate,1]{label=\alph*),ref=\thetheorem(\alph*)}

\makeatletter
\def\paragraph{\vskip0.5ex\@startsection{paragraph}{4}%
  \z@\z@{-\fontdimen2\font}%
  {\normalfont\bfseries}}
\def\subparagraph{\vskip0.5ex\@startsection{subparagraph}{5}%
  \z@\z@{-\fontdimen2\font}%
  {\normalfont\bfseries}}
\makeatother

\begin{document}
\maketitle

\begin{abstract}
The Tur\'an problem asks for the largest number of edges in an $n$-vertex graph
not containing a fixed forbidden subgraph $F$. We construct a new family
of graphs not containing $K_{s,t}$, for $t= C^s$, with $\Omega(n^{2-1/s})$ edges matching the upper bound
of K\"ov\'ari, S\'os and Tur\'an. 
\end{abstract}

\section{Introduction}
\paragraph{The Tur\'an problem.}
Let $F$ be a fixed graph. The Tur\'an problem asks for the value of $\ex(n,F)$, the largest
number of edges in an $n$-vertex graph not containing a copy of $F$ as a subgraph. The classic
theorem of Erd\H{o}s and Stone \cite{erdos_stone} gives an asymptotic for $\ex(n,F)$ when $F$ is
not bipartite.

For bipartite $F$, much less is known. Even the simplest case when $F$ is a complete bipartite graph $K_{s,t}$
is open. Specifically, K\"ov\'ari, S\'os and Tur\'an \cite{kst} proved that 
\[
  \ex(n,K_{s,t})=O_{s,t}(n^{2-1/s}).
\]
Obviously, we may reverse the roles of $s$ and $t$ to obtain $\ex(n,K_{s,t})=O_{s,t}(n^{2-1/t})$, which is superior if $t<s$.
So, from now on we discuss only the case $t\geq s$.
Though the implicit constant in the big-Oh notation has been improved by F\"uredi \cite{furedi_zarankiewicz_upper}, the K\"ov\'ari--S\'os--Tur\'an bound remains
the only upper bound on $\ex(n,K_{s,t})$. Many researchers conjecture that the K\"ov\'ari--S\'os--Tur\'an bound
is tight (e.g., \cite[p.~52]{kst}, \cite[p.~6]{erdos_conj}, \cite[p.~257]{furedi_turansurvey}, \cite[Conjecture~2.24]{furedi_simonovits_survey}).
However, apart from the numerous results for $s=2$ and $s=3$ (see \cite[Section~3]{furedi_simonovits_survey} for a survey),
there are only two constructions attaining the K\"ov\'ari--S\'os--Tur\'an bound for general $s>3$. The first is due to Alon, R\'{o}nyai and Szab\'{o} \cite{alon_ronyai_szabo}
who, improving on the previous construction by Koll\'ar, R\'onyai and Szab\'o \cite{kollar_ronyai_szabo}, showed that
\begin{equation}\label{eq:sharpkst}
  \ex(n,K_{s,t})=\Omega_s(n^{2-1/s})\qquad\text{if }t>(s-1)!.
\end{equation}
The construction is a clever use of norms over finite fields. 
The second, more recent class of constructions originating from \cite{blagojevic_bukh_karasev} uses random varieties.
It has hitherto provided inferior dependence of $t$ on~$s$. For example, \cite{blagojevic_bukh_karasev} obtains \eqref{eq:sharpkst}
only for $t\geq s^{4s}$. The advantage of these constructions is their flexibility, see \cite{conlon_theta,bukh_conlon,ma_co,he_tait,bukh_tait,xu_zhang_ge} for some of their variations and applications.

In this work, we use a novel version of the random algebraic method to construct graphs that match the
K\"ov\'ari--S\'os--Tur\'an bound for $t$ that is only exponential in~$s$.
\begin{theorem}\label{thm:turan}
Let $s\geq 2$. Then
\[
  \ex(n,K_{s,t})=\Omega_s(n^{2-1/s})\qquad\text{if }t>9^s\cdot s^{4s^{2/3}}.
\]
\end{theorem}

\paragraph{The Zarankiewicz problem.}
Closely related to the Tur\'an problem for $K_{s,t}$-free graphs is the problem of Zarankiewicz \cite{zarankiewicz}.
It is the asymmetric version of the Tur\'an problem. It is well-known that in the study of the growth rate of $\ex(n,F)$
we may assume that the $n$-vertex graph is bipartite. To distinguish the two parts of a bipartite graph, we shall
call them the \emph{left} and \emph{right}. A copy of $K_{s,t}$ in a bipartite graph $G$ can be situated in two ways:
either the $s$ vertices are in the left part, or the $s$ vertices are in the right part. In the Zarankiewicz problem, we forbid only the former case.
So, we say that $G$ is a \emph{sided graph} if it is bipartite with distinguished left and right parts.
The Zarankiewicz problem then asks for the estimate on the number of edges in a sided graph
not containing $K_{s,t}$, which is regarded as a sided graph with $s$ vertices on the left and $t$ vertices on the right.

An important consequence of making the graph bipartite with distinguished parts is that the two parts can (possibly) be 
of very unequal size. This often occurs in applications (see e.g.~\cite{alon_mellinger_mubayi_verstraete,walsh}). With this in mind, define $z(m,n;s,t)$ as the largest number
of edges in a sided graph with $m$ vertices on the left and $n$ vertices on the right that contains no sided~$K_{s,t}$.
The bound of K\"ov\'ari, S\'os and Tur\'an for the Zarankiewicz problem takes the form
\[
  z(m,n;s,t)=O_{s,t}(mn^{1-1/s}).
\]
In the symmetric case when $m=n$, the best known constructions for Zarankiewicz problem are the same
as the best bipartite constructions for the Tur\'an problem. This is not so for our approach: we
are able to take the advantage of the fact that only one orientation of $K_{s,t}$ is forbidden
to obtain a lower bound on $z(n,n;s,t)$ that is superior to the corresponding bound for $\ex(n,K_{s,t})$
in \Cref{thm:turan}.
\begin{theorem}\label{thm:zar}\begin{enumerate}
\item \label{thm:zarsmall}
  Suppose $s,t,m,n,k\geq 3$ are integers that satisfy the inequalities $\log_n m\leq\nobreak \frac{s^{k-2}}{k!\log^{2k} s}$ and $k\leq \frac{s}{2\log^3 s}$. Then
  \[
    z(m,n;s,t)=\Omega_{s}(mn^{1-1/s})\qquad\text{if }t>k^s\cdot e^{2s/\log s}.
  \]
In particular, $z(n,n;s,t)=\Omega_s(n^{2-1/s})$ for $t>3^{s+o(s)}$.

\item \label{thm:zarlarge}
For each $s\geq 3$ there is a constant $c_s>0$ such that if the inequality $\log_n m\leq c_st^{\frac{1}{1+2\log s}}$ holds, then
\[
  z(m,n;s,t)=\Omega_{s}(mn^{1-1/s}).
\]
\end{enumerate}
\end{theorem}
Part (b) is an improvement on the result of Conlon \cite{conlon_remarks}, who proved the $\Omega_s(mn^{1-1/s})$ bound under
the condition $\log_n m\leq c_s' t^{1/(s-1)}$. In the same article, Conlon asked if the bound holds for $\log_n m\leq t/s$,
which would be tight if true.

\paragraph{Main proof idea, in a nutshell.} The key novelty in our construction is that it is `bumpy'. All the 
previous algebraic constructions, whether random or not, were flat: the vertex of the graph was 
always an affine or a projective space over $\Fq$, which was occasionally mildly mutilated by having a lower-dimensional
subset removed, or stitched to itself via a quotient operation. These are finite field analogues of $\mathbb{R}^n$ with a flat metric. 
In contrast, the vertex set in our construction is a solution set to a family of random polynomial equations,
and cannot be flattened with any change of coordinates.

\paragraph{Proof ideas, in more detail.} To explain our construction, we recall the important ingredients
in the previous random algebraic constructions of $K_{s,t}$-free graphs. 
The key is an estimate on the probability that the variety $\V(f_1,\dotsc,f_s)$ cut out by $s$ random polynomials
contains many points. This relies on two inputs. The first is Bezout's theorem, which is used to conclude that
\emph{if} $\V(f_1,\dotsc,f_s)$ is zero-dimensional, then it contains at most $\prod \deg f_i$ many points.
The second input is a bound on the probability that $\V(f_1,\dotsc,f_s)$ is of codimension less than~$s$.

To prove that bound, the original paper \cite{blagojevic_bukh_karasev} uses Hilbert functions
(though it did not use probabilistic language). However, almost all the subsequent
works use the approach from \cite{bukh_random} relying on the Lang--Weil bounds, the only exception
being an elegant argument in \cite{conlon_remarks} which can be described as an implicit use of Hilbert functions.

In this paper, we go back to the explicit use of Hilbert functions. We show
that the codimension of $\V(f_1,\dotsc,f_s)$ is extremely likely equal to $s$, 
\emph{unless} $s$ is close to the dimension of the ambient space. That is precisely
the situation in the previous works, where the graph's vertex set was $s$-dimensional set $\Fq^s$.
To bypass this obstacle we construct the vertex set in two steps: At the start we use
affine space $\Fq^{s+r}$ of slightly larger dimension. This way the varieties of the form $\V(f_1,\dotsc,f_s)$ that we obtain 
have dimension $r$ with very high probability.
We then shrink the vertex set to a random subvariety of codimension $r$, thereby cutting all the varieties
of the form $\V(f_1,\dotsc,f_s)$ at once.

Our second innovation concerns the need to control $\prod \deg f_i$ in the Bezout's bound. In the construction of Tur\'an
graphs, the random polynomials $f_1(y),\dotsc,f_s(y)$ arise as specializations of a single polynomial $g(x,y)$ in
two sets of variables. A simple way to ensure that the random polynomials $f_1(y),\dotsc,f_s(y)$ are mutually independent
is to make $g$ a random polynomial of degree at least~$s$. That makes $\prod \deg f_i$ grow like~$s^s$. To circumvent
this, we further replace the $(s+r)$-dimensional affine space by a variety $U$ that has the property that
every specialization of a random polynomial $g$ of bounded degree to any $s$ points of $U$ yields $s$ mutually
independent polynomials. We call these \emph{$m$-independent varieties}. The construction of $m$-independent varieties occupies the 
bulk of the paper (\Cref{sec:constindep,sec:vdt}).\medskip

Finally, we want to highlight an auxiliary contribution of this work that is of independent interest. The construction of $m$-independent
varieties depends on a bound on the number of minimal linear dependencies
between $m$'th powers of linear forms, i.e., $0=\sum \alpha_i\ell_i(x)^m$ where $\ell_i(x)=c_{i,0}x_0+\dotsb+c_{i,b}x_b$. Here, `minimal' means
that no proper subset of the $m$'th powers of these linear forms is linearly dependent. Representation of polynomials by
sums of $m$'th powers of linear forms has been much studied, motivated primarily by the Waring problem. In particular,
we adapt the argument which is implicit in the work of Bia{\l}ynicki\nobreakdash-Birula and Schinzel \cite{bialynicki_schinzel} to our purposes.
Though the argument in \cite{bialynicki_schinzel} suffices to obtain an exponential bound in \Cref{thm:turan}, we go beyond and obtain a stronger bound that yields
the smaller exponent base of~$9$ in \Cref{thm:turan}.

\paragraph{Paper organization.} We begin by collecting the algebraic tools we require in \Cref{sec:tools}.
The concept of an $m$-independent set, which is central to the proof of \Cref{thm:turan}, is introduced in
\Cref{sec:defindep}. The $m$-independent varieties used in the proof of \Cref{thm:turan} are constructed in \Cref{sec:constindep,sec:vdt}.
Finally, in \Cref{sec:turan,sec:zar} we prove \Cref{thm:turan,thm:zar} respectively.

\paragraph{Acknowledgments.} I am thankful to Chris Cox for extensive feedback on an earlier version
of this paper. I am grateful to Jacob Tsimerman for discussions on numerous topics related to this paper, and especially for motivating me to prove
\Cref{lem:mdependentdim} in its current form, and to Anamay Tengse, Mrinal Kumar, Ramprasad Saptharishi for spotting a serious error in
the proof of \Cref{lem:mdependentdim} in the previous version of this paper. Finally, I am grateful to the anonymous referees for a number of
constructive comments.

\section{Algebraic tools}\label{sec:tools}
To make this paper maximally accessible, we tried to keep the use of algebra to the minimum. In particular, we use
counting arguments even when similar algebraic arguments could have provided slightly superior numeric constants. 
Despite this, we require basic familiarity with algebraic geometry on the level of the first chapter of Shafarevich's book \cite{shafarevich}. 
We collect the other algebraic tools in this section.

\paragraph{Varieties and their \texorpdfstring{$\Fq$-points}{F\textunderscore q-points}.}
The integer $q$ will denote a prime power. We shall work exclusively with fields $\Fq$ and $\Fqbar$.
All varieties in this paper are quasi-projective over the field~$\Fqbar$.
We write $\V(f_1,\dotsc,f_t)$ for the projective variety cut out by homogeneous polynomials
$f_1,\dotsc,f_t$.
We denote the vector space of homogeneous polynomials of degree $m$ in $b+1$ variables with coefficients in $\Fq$ by $\poly[m]{\Fq}{x}{b}$.
We also work with products of projective spaces $\P^a\times \P^b$. The set of
bihomogeneous polynomials of bidegrees $(m,m')$ on $\P^a\times \P^b$ is denoted by
$\poly[m]{\Fq}{x}{a}\otimes\nobreak \poly[m']{\Fq}{y}{b}$.

We write monomials using the multiindex notation: for a multiindex $\beta=(\beta_0,\beta_1,\dotsc,\beta_b)\in \Znonneg^{b+1}$,
the notation $x^{\beta}$ stands for the monomial $x_0^{\beta_0}x_1^{\beta_1}\dotsb x_b^{\beta_b}$. A general homogeneous polynomial
of degree~$m$ is thus written as $\sum_{\abs{\beta}=m} c_{\beta}x^{\beta}$.

The graphs we shall construct in the proofs of \Cref{thm:turan,thm:zar} will consist of the
$\Fq$-points of certain varieties. We denote the $\Fq$-points of a variety $V\subseteq \P^b$ by $V(\Fq)$ or (if the
variety $V$ is a complicated expression) by $V\cap \P^b(\Fq)$. We shall use the following bounds
on the number of~$\Fq$-points.

\begin{lemma}[Weakening of {\cite[Corollary 3.3]{couvreur}}]\label{zippel}
Suppose $V\subseteq \P^b$ is any $k$\nobreakdash-dimensional variety of degree~$d$. Then $\abs{V(\Fq)}\leq d\abs{\P^k(\Fq)}$.
\end{lemma}

\begin{lemma}\label{lem:expected}
  \begin{enumerate}
  \item \label{part:expectproj}
  Let $m_1,\dotsc,m_r$ be positive integers, and let $Y\subseteq \P^b(\Fq)$ be a non-empty set. Suppose that $g_1,\dotsc,g_r\in \poly{\Fq}{x}{b}$ are random homogeneous polynomials of
  degrees $\deg g_i=m_i$. Then
  \begin{equation}\label{eq:expected}
    \Pr\left[\abs{Y\cap \V(g_1,\dotsc,g_r)}\leq \frac{\abs{Y}}{2q^r}\right]\leq \frac{4q^r}{\abs{Y}}. 
  \end{equation}
  \item \label{part:expectedbi} The same holds for bihomogeneous polynomials, i.e., if $Y\subseteq \P^a(\Fq)\times \P^b(\Fq)$ is a non-empty set and $g_1,\dotsc,g_r$ are random bihomogeneous polynomials of bidegrees
    $\deg g_i=(m_i,m_i')$ with $m_i,m_i'\geq \nobreak 1$, then \eqref{eq:expected} holds.
  \end{enumerate}
\end{lemma}
\begin{proof}
  For a point $y\in Y$, let $R_y$ be the indicator random variable of the event $y\in \V(g_1,\dotsc,g_r)$.
  Let $y,y'\in Y$ be two distinct points. We claim that $\E[R_y]=\E[R_{y'}]=1/q^r$ and that the random variables $R_y$ and $R_{y'}$
  are independent. To see this in the \hyperref[part:expectproj]{case (a)}, apply a change of coordinates so that 
  $y=[1:0:0:\dotsc:0]$ and $y'=[0:1:0:\dotsc:0]$. The polynomial $g_i$ vanishes at $y$ if and only if the coefficient
  of $x_0^m$ vanishes. Similarly, $g_i(y')=0$ if and only if the coefficient of $x_1^m$ vanishes.
  In the \hyperref[part:expectedbi]{case (b)}, write $y=(y_a^{},y_b^{})$ and $y'=(y_a',y_b')$ with $y_a^{},y_a'\in \P^a(\Fq)$ and $y_b^{},y_b'\in \P^b(\Fq)$.
  Since $y\neq y'$ we may assume that $y_b^{}\neq y_b'$ (by swapping the roles of $\P^a$ and $\P^b$ if necessary). We can then
  change the coordinates on $\P^b$ in the same way as in the \hyperref[part:expectproj]{case (a)}, and observe that the vanishing of $g_i$ at $y$ and $y'$ depends on disjoint sets
  of coefficients. This proves the claim.\smallskip

  Let $R=\sum_{y\in Y} R_y$. From the pairwise independence of the $R_y$'s, it follows that 
  \[
    \Var[R]=\sum_{y\in Y} \Var[R_y]=\tfrac{1}{q^r}(1-\tfrac{1}{q^r})\abs{Y}\leq \abs{Y}/q^r.
  \]
  Since $\E[R]=\sum_y \E[R_y]=\abs{Y}/q^r$, Chebyshev's inequality then implies that
  \[
    \Pr\bigl[\abs{R-\abs{Y}/q^r}\geq \lambda\sqrt{\abs{Y}/q^r} \bigr]\leq \lambda^{-2},
  \]
  and the lemma follows upon taking $\lambda=\sqrt{\abs{Y}/4q^r}$.
\end{proof}

\paragraph{B\'ezout's inequality.} For a reducible variety $V\subseteq \P^b$, define the \emph{total degree} $\deg(V)$ to be the sum
of the degrees of the irreducible components of~$V$.
\begin{lemma}[B\'ezout's inequality, {\cite[p.~228, Example 12.3.1]{fulton}}]
  Let $V$ and $W$ be two varieties in $\P^b$. Then
  \[
    \deg(V\cap W)\leq \deg(V)\deg(W).
  \]
\end{lemma}

\paragraph{Hilbert functions.} For a homogeneous ideal $I$, the Hilbert function $H_I(m)$ is defined as the codimension
of $I_m\eqdef I\cap \poly[m]{\Fqbar}{x}{b}$ in $\poly[m]{\Fqbar}{x}{b}$. Equivalently, if $I=I(V)$ is the homogeneous ideal of
polynomials vanishing on a variety $V$, then $H_I(m)$ is the dimension of the subspace of functions on $V$ induced by all homogeneous polynomials of degree~$m$.

We shall use the following bound on the Hilbert function.
\begin{lemma}\label{lem:hilbert}
Let $V$ be a variety of dimension $k$ in $\P^b$. Then its Hilbert function satisfies
\[
  H_{I(V)}(m) \geq \binom{m+k}{m}.
\]
\end{lemma}
Though more general bounds are known (e.g.~\cite[Theorem~2.4]{sombra}), we give a proof for completeness.
\begin{proof}
  By \cite[Theorem~1.15 and Corollary~1.6]{shafarevich} there is a linear projection $\pi\colon V\to \P^k$, which is finite and surjective.
  Then the pullback of a homogeneous polynomial of degree $m$ on~$\P^k$ is a homogeneous
  polynomial of degree $m$ on~$V$. The result follows since the space of homogeneous degree-$m$ polynomials on $\P^k$ is of dimension $\binom{m+k}{m}$,
  and the pullback map has trivial kernel.
\end{proof}
We use Hilbert functions to bound the probability that a random polynomial vanishes on a given variety.
\begin{lemma}\label{use:hilbert}
  Let $V$ be any variety in~$\P^b$. Let $g\in \poly[m]{\Fq}{x}{b}$ be a random homogeneous degree-$m$ polynomial.
  Then the probability that $g$ vanishes on $V$ is
  \[
    \Pr[g_{|V}=0]\leq q^{-H_{I(V)}(m)}.
  \]
\end{lemma}
\begin{proof}
  By the definition of the Hilbert function, $H_{I(V)}(m)$ is the codimension (over $\Fqbar$)
  of $I(V)_m$ in $\Fqbar[x_0,\dotsc,x_b]_m$. This implies that the codimension of $I(V)_m \cap \Fq[x_0,\dotsc,x_b]$
  in $\Fq[x_0,\dotsc,x_b]_m$ is at least $H_{I(V)}(m)$. Since $\Pr[g_{|V}=0]=\Pr[g\in I(V)_m \cap \Fq[x_0,\dotsc,x_b]]$,
  the lemma follows.
\end{proof}

\section{Definition and uses of \texorpdfstring{$m$-independence}{m-independence}}\label{sec:defindep}
If $t$ is small compared to $b$, then the Hilbert function
of $t$ generic points in $\P^b$ is $H_{I(\{p_1,\dotsc,p_t\})}(m)=t$ for $m\geq 1$.
As we shall see shortly, the point sets satisfying $H_{I(\{p_1,\dotsc,p_t\})}(m)=t$
behave random-like with respect to the degree-$m$ polynomials. Since this property
will be key to our construction of Tur\'an graphs, we focus on the sets
lacking this pleasant property.

\begin{definition}\label{def:mdep}
  We say that points $p_1,\dotsc,p_t\in \P^b$ are \emph{$m$-dependent}
  if $H_{I(\{p_1,\dotsc,p_t\})}(m)<t$. We furthermore say that they are
  \emph{minimally $m$\nobreakdash-dependent} if no proper subset of these points is $m$\nobreakdash-dependent.
\end{definition}
\begin{definition}\label{def:mwisedep}
A set $X\subset \P^b$ is called \emph{$s$-wise $m$-independent} if no $s$ distinct points of $V$ are $m$\nobreakdash-dependent.
\end{definition}

%

Though we define $m$-dependence via Hilbert functions, there are two other ways to think about the concept that will be useful.
The first way is to think of a homogeneous degree-$m$ polynomial $f$ on $\P^b$ as a linear form in its coefficients.
Specialization of $f$ to $f(p)$ gives distinct linear forms for distinct $p\in \P^b$. It is easy
to see from the definition of $m$-dependence that the points $p_1,\dotsc,p_t$
are $m$-dependent if and only if $f(p_1),\dotsc,f(p_t)$ are linearly dependent as linear forms in $\binom{b+m}{m}$ variables.

The second way to understand $m$-dependence is via projective space duality, which we think of as an identification between points of $\P^b$ and linear forms.
Namely, to each point $p\in \P^b$ we associate the linear form $\ell$ in $b+1$ variables defined by $\ell(x)\eqdef \langle x,p\rangle$.
Because points in the projective space are defined only up to a multiplication by a non-zero scalar, this linear form too is
defined up to a multiplication by a non-zero scalar.
\begin{lemma}\label{eq:hilbequiv}
Let $p_1,\dotsc,p_t\in \P^b$ be any $t$ points, and let $\ell_1(x),\dotsc,\ell_t(x)$ be the linear forms associated to these points.
The following are equivalent:
\begin{enumerate}
\item points $p_1,\dotsc,p_t$ are $m$-dependent,
\item there is a linear relation of the form
\begin{equation}\label{eq:linprodrel}
  c_1 \prod_{j=1}^m \ell_1(x_j)+\dotsb+c_t \prod_{j=1}^m \ell_t(x_j)=0.
\end{equation}
\end{enumerate}
Furthermore, if the characteristic of $\Fq$ is at least $m$, then the following condition is also equivalent to
the two above:
\begin{enumerate}[resume]
\item there is a linear relation of the form
\begin{equation}\label{eq:linpowrel}
  c_1 \ell_1(x)^m+\dotsb+c_t \ell_t(x)^m=0.
\end{equation}
\end{enumerate}
\end{lemma}
\begin{proof}[Proof of (a)$\iff$(b).]\let\qed\relax
Let $k_1,\dotsc,k_m\in\{0,1,\dotsc,b\}$ be arbitrary, and let $\beta=(\beta_0,\beta_1,\dotsc,\beta_b)$ be the tally vector,
i.e., $\beta_i$ is the number of elements among $k_1,\dotsc,k_m$ that are equal to~$i$. Then
the coefficient of $x_{1,k_1}\dotsb x_{m,k_m}$ in $\prod_j \ell_i(x_j)$ is equal to $p_i^{\beta}$.
Therefore, the linear relation in \eqref{eq:linprodrel} holds if and only if 
$c_1 p_1^{\beta}+\dotsb+c_t p_t^{\beta}=0$ for every $\beta\in \Znonneg$ satisfying $\abs{\beta}=m$.
As discussed in the paragraph following \Cref{def:mdep,def:mwisedep}, the latter condition is equivalent to 
points $p_1,\dotsc,p_t$ being $m$-dependent.
\end{proof}
\begin{proof}[Proof of (b)$\implies$(c).]\let\qed\relax Trivial, by setting $x_1=\dotsb=x_m$.\end{proof}
\begin{proof}[Proof of (c)$\implies$(b).] Set $x=x_1+\dotsb+x_m$ in \eqref{eq:linpowrel}, and expand. 
Let $k_1,\dotsc,k_m$ be arbitrary, and define $\beta$ as in the proof of (a)$\iff$(b).
The coefficient of $x_{1,k_1}\dotsb x_{m,k_m}$ in the expansion is $\binom{m}{\beta_0,\beta_1,\dotsc,\beta_b}$
times larger than the coefficient of the same monomial in \eqref{eq:linprodrel}. Since $\binom{m}{\beta_0,\beta_1,\dotsc,\beta_b}\neq 0$
follows from the assumption on the characteristic, the desired implication follows. 
\end{proof}

If points $p_1,\dotsc,p_t$ are not $m$-dependent, we say that they are \emph{$m$-independent}.

The usefulness of these definitions comes from the combination of two simple observations:
\begin{proposition}\label{prop:indep}
Suppose that $v_1,\dotsc,v_s\in \P^a(\Fq)$ are $m$-independent points. Pick a random polynomial $g\in \nobreak\poly[m]{\Fq}{x}{a}\otimes \poly[m']{\Fq}{y}{b}$ uniformly among all bihomogeneous polynomials 
of bidegree $(m,m')$ on $\P^a\times\P^b$. 
Then the $s$ random polynomials $g(v_1,y),\dotsc,g(v_s,y)\in \poly[m]{\Fq}{y}{b}$ are mutually independent.
\end{proposition}
\begin{proof}
Write $g$ as $g(x,y)=\sum_{\abs{\beta}=m'} y^{\beta}g_{\beta}(x)$. Since the coefficient of $y^{\beta}$ in $g(v_i,y)$
is $g_{\beta}(v_i)$, it suffices to show that $g_{\beta}(v_1),\dotsc,g_{\beta}(v_s)$ are independent
for every $\beta$. 

Think of $g_{\beta}(v_i)$ as a linear function of the coefficients of $g_{\beta}$.
By the alternative definition of $m$-independence discussed above, the linear functions $g_{\beta}(v_1),\dotsc,g_{\beta}(v_s)$
are linearly independent. Since the coefficients of $g_{\beta}$ are chosen independently and
uniformly from $\Fq$, this implies that the random variables $g_{\beta}(v_1),\dotsc,g_{\beta}(v_s)$ are independent.
\end{proof}

Define the function
\[
  M_k(t)\eqdef \min\Bigl\{ m :\binom{m+k}{k}\geq t\Bigr\}.
\]
Recall that a variety $W$ is said to be of \emph{pure dimension} $k$ if all of its irreducible components
are of dimension $k$.
\begin{proposition}\label{prop:randomcut}
Let $0\leq s\leq k\leq b$ be integers, and
let $W\subseteq \P^b$ be a variety of pure dimension~$k$ and degree at most~$D$.
\begin{enumerate}
\item \label{cut:same}
If $h_1,\dotsc,h_s$ are independent random homogeneous polynomials of degree $m$ in $b+1$ variables with $\Fq$-coefficients,
then 
\[
  \Pr\bigl[\dim(W\cap \V(h_1,\dotsc,h_s))>k-s\bigr]\leq Cq^{-\binom{k-s+1+m}{m}},
\]
where the constant $C=C(k,m,D)>0$ depends only on $k,m$ and $D$.
\item \label{cut:diff}
Let $T$ and $\delta_1,\dotsc,\delta_s$ be positive integers satisfying $\delta_i\geq M_{k-i+1}(T)$.
If $h_1,\dotsc,h_s$ are independent random homogeneous polynomials of degrees $\deg h_i=\delta_i$ in $b+1$ variables with $\Fq$-coefficients,
then 
\[
  \Pr\bigl[\dim(W\cap \V(h_1,\dotsc,h_s))>k-s\bigr]\leq Cq^{-T},
\]
where the constant $C=C(T,s,D,\delta_1,\dotsc,\delta_s)>0$ depends only on $T,s$,$D$ and on the polynomial degrees~$\delta_1,\dotsc,\delta_s$.
\end{enumerate}
\end{proposition}
\begin{proof}
\hyperref[cut:same]{Part (a)} is the special case of \hyperref[cut:diff]{part (b)} with $T=\binom{k-s+1+m}{m}$ and $\delta_1=\dotsb=\delta_s=m$. So, it suffices to prove \hyperref[cut:diff]{part (b)}.

The proof is by induction on~$s$. Let $\mathcal{U}$ be the set of all irreducible components
of $W\cap \V(h_1,\dotsc,h_{s-1})$. Since $\V(h_1,\dotsc,h_{s-1})$ is cut out by $s-1$ polynomials and $W$ is of pure dimension~$k$, it follows that $\dim(U)\geq k-s+1$ for
each $U\in\mathcal{U}$ \cite[Corollary 1.14]{shafarevich}. So, from \Cref{lem:hilbert,use:hilbert} we see that, for any fixed $U\in\mathcal{U}$,
\[
  \Pr[h_s\text{ vanishes on }U]\leq q^{-\binom{k-s+1+\delta_s}{k-s+1}}\leq q^{-\binom{k-s+1+M_{k-s+1}(T)}{k-s+1}}\leq q^{-T}.
\]
B\'ezout's inequality tells us that $\abs{\mathcal{U}}\leq D\prod_{i=1}^{s-1} \delta_i$, and hence
the probability that $h_s$ vanishes on some component of $W\cap V(h_1,\dotsc,h_{s-1})$ is at most $D\prod_{i=1}^{s-1} \delta_i\cdot q^{-T}$.
By the induction hypothesis, the variety $W\cap \V(h_1,\dotsc,h_{s-1})$ is of dimension exceeding $k-s+1$
with probability at most $C(T,s-\nobreak 1,D,\delta_1,\dotsc,\delta_{s-1})q^{-T}$, and so the probability that the variety $W\cap \V(h_1,\dotsc,h_{s-1},h_s)$
is of dimension exceeding $k-s$ is at most
\[
  D\prod_{i=1}^{s-1} \delta_i\cdot q^{-T}+C(t,s-1,D,\delta_1,\dotsc,\delta_{s-1})q^{-T},
\]
which is at most $C(T,s,D,\delta_1,\dotsc,\delta_s)q^{-T}$ for a suitable $C(\dotsb)$.
\end{proof}

Combining these two observations we obtain the following handy result.
\begin{lemma}\label{lem:handy}
Let $0\leq s\leq k\leq b$ be integers.
Suppose that $W\subseteq \P^b$ is a variety of pure dimension~$k$, and
$X\subset \P^a(\Fq)$ is an $s$-wise $m$\nobreakdash-independent set of size
$\abs{X}\leq c' q^{\frac{1}{s}\binom{k-s+1+m}{m}}$ where $c'=c'(m,s,\deg W)>0$ is sufficiently small.
Let $g\in \poly[m]{\Fq}{x}{a}\otimes \poly[m]{\Fq}{y}{b}$ be a random bihomogeneous polynomial 
of bidegree $(m,m)$ on $\P^a\times \P^b$. Then the following holds with probability at least $\tfrac{4}{5}$:
For every $s$ distinct points $v_1,\dotsc,v_s\in X$ the variety
\[
  \{w\in W : g(v_1,w)=\dotsb=g(v_s,w)=0\}
\]
is of dimension $k-s$.
\end{lemma}
\begin{proof}
  By \Cref{prop:indep}, polynomials $g(v_1,w),\dotsc,g(v_s,w)$ are mutually independent, for every $s$ distinct points $v_1,\dotsc,v_s\in X$.
  By \Cref{cut:same} and the union bound over all $(v_1,\dotsc,v_s)\in X^s$, the probability that $g$ does not satisfy the lemma
  is at most
  $\abs{X}^s\cdot C(m,s,\deg W) q^{-\binom{k-s+1+m}{m}}\leq (c')^s C(m,s,\deg W)$. So, choosing $c'$ small enough works.
\end{proof}

\section{Construction of \texorpdfstring{$m$-independent}{m-independent} varieties}\label{sec:constindep}
For positive integers $b,m,t$, let 
\begin{align*}
  \VD_t(b,m)&\eqdef \{(p_1,\dotsc,p_t)\in (\P^b)^t : p_1,\dotsc,p_t\text{ are minimally $m$-dependent} \},\\
  \vd_t(b,m)&\eqdef \dim \VD_t(b,m).
\end{align*}
Observe that $\VD_t(b,m)$ is a variety. Indeed, for a subset $I\subseteq [t]\eqdef \{1,2,\dotsc,t\}$, define
\[
  \widetilde{\VD}_I(b,m)\eqdef \{(p_1,\dotsc,p_t)\in (\P^b)^t : (p_i)_{i\in I}\text{ are $m$-dependent}\}.
\]
The set $\widetilde{\VD}_I(b,m)$ is a projective variety because of the equivalence of $m$-dependence of points and linear
dependence of the linear $m$-th powers of respective linear forms, which we discussed above.
From this, it follows that $\VD_t(b,m)=\widetilde{\VD}_{[t]}(b,m)\setminus \bigcup_{I\subsetneq [t]} \widetilde{\VD}_I(b,m)$ is a difference between two projective varieties,
and so itself is a (quasiprojective) variety.

We shall upper bound the functions $\vd_t$ in the next \namecref{sec:vdt}. But first we show how to use these bounds to construct $m$-independent varieties.
\begin{lemma}\label{lem:fqbarindep}
  Suppose that $b,m,Z,s\geq 1$ are integers satisfying
  \begin{equation}\label{eq:Zcond}
    Z>\frac{1}{t-1}\vd_t(b,m)\qquad\text{ for all }t=2,3,\dotsc,s.
  \end{equation}
  Let $f_1,\dotsc,f_{Z}$ be generic degree-$m$ homogeneous polynomials in $b+1$ variables.
  Then the variety $\V(f_1,\dotsc,f_{Z})$ is $s$-wise $m$-independent.
\end{lemma}
As we shall show in \Cref{part:mdep_empty}, $\VD_t(b,m)=\emptyset$ for $t\leq m+1$. Hence, the condition \eqref{eq:Zcond} is non-vacuous only
for $t=m+2,m+3,\dotsc,s$.
\begin{proof}[Proof of \Cref{lem:fqbarindep}]
  It suffices to show that $\V(f_1,\dotsc,f_Z)$ contains no set of $t$ minimally $m$\nobreakdash-dependent points
  for every $t=2,3,\dotsc, s$. Fix $t$ in this range.

  Note that whenever points $p_1,\dotsc,p_t\in \P^b$ are minimally $m$\nobreakdash-dependent, their Hilbert function satisfies $H_{I(\{p_1,\dotsc,p_t\})}(m)= t-1$.
  Indeed, the inequality $H_{I(\{p_1,\dotsc,p_t\})}(m)\leq t-1$ follows from the definition of $m$\nobreakdash-dependence, and $H_{I(\{p_1,\dotsc,p_t\})}(m)\geq t-1$ follows from minimality.
  This means that the vector space $I(\{p_1,\dotsc,p_t\})_{m}$ is of codimension $t-1$ in $\Fqbar[x_0,\dotsc,x_b]_{m}$,
  and hence
  $I(\{p_1,\dotsc,p_t\})_{m}^Z$
  is of codimension $Z(t-1)$ in $(\Fqbar[x_0,\dotsc,x_b]_{m})^Z$.
  
  Since $Z(t-1)>\vd_t(b,m)$, we can use the algebraic version of the union bound to deduce that
  for generic degree-$m$ homogeneous polynomials $f_1,\dotsc,f_Z$, the variety $\V(f_1,\dotsc,f_Z)$ does not contain any minimally $m$-dependent set $\{p_1,\dotsc,p_t\}$.
  Indeed, write $F\eqdef\Fqbar[x_0,\dotsc,x_b]_{m}^Z$ and regard it as a variety with polynomials' coefficients as indeterminants.
  Define the variety
  \[
    V\eqdef \{(f_1,\dotsc,f_Z,p_1,\dotsc,p_t)\in F\times \VD_t(b,m) : f_i(p_j)=0\text{ for all }i,j\}.
  \]
  Consider the projection of $V$ onto~$F$. Our claim is that the fiber of a generic $\vec{f}\in F$ is empty. If this
  is not so, then $\dim V\geq \dim F$. However, for the projection of~$V$ onto the $\VD_t(b,m)$ factor, every fiber
  is of codimension $Z(t-1)$, and so $\dim V\leq \dim \VD_t(b,m)+(\dim F-Z(t-1))<\dim F$, which is a contradiction,
  proving our claim.
\end{proof}
\begin{lemma}\label{lem:fqindep}
  Suppose that $b,m,Z,s$ are integers satisfying the condition \eqref{eq:Zcond} and $b-Z\geq 1$, and assume that 
  $q\geq q_0(b,m,Z)$ is sufficiently large in terms of $b,m,Z$.
  Then there exist degree-$m$ polynomials $f_1,\dotsc,f_Z\in \poly[m]{\Fq}{x}{b}$
  such that the variety $\V(f_1,\dotsc,f_{Z})$ is $s$-wise $m$-independent, is of pure dimension $b-Z$, and contains
  at least $\tfrac{1}{2}q^{b-Z}$ many $\Fq$-points.
\end{lemma}
\begin{proof}
We shall choose each $f_1,\dotsc,f_Z$ uniformly at random from $\poly[m]{\Fq}{x}{b}$.
  Let 
\begin{align*}
  B&\eqdef \{(f_1,\dotsc,f_Z) : \V(f_1,\dotsc,f_{Z})\text{ is not }s\text{-wise }m\text{-independent}\}\\
   &\subseteq \poly[m]{\Fqbar}{x}{b}^Z.\\
\intertext{The set $B$ is a variety, for we can think of it as the image of the projection of the bigger variety}
  B'&\eqdef\{(f_1,\dotsc,f_Z, p_1,\dotsc,p_t) : f_i(p_j)=0\text{ for all }i,j\}\\
    &\subseteq \poly[m]{\Fqbar}{x}{b}^Z\times \VD_t(b,m)
\end{align*}
onto the first factor.

Furthermore, recall that $\VD_t(b,m)$ is defined by the linear dependence of the $m$-th powers of
the linear forms associated to the points $p_i$. Since the number of linear forms and the numbers of variables therein
do not depend on $q$, using B\'ezout's inequality we may obtain an upper bound on the degree $B'$ (and hence on the degree of $B$)
that is independent of $q$ (but depends on $b,m$ and $Z$). 

By \Cref{lem:fqbarindep}, $B$ is of codimension at least $1$ in $\poly[m]{\Fqbar}{x}{b}^Z$.
By \Cref{zippel}, a random element $\poly[m]{\Fq}{x}{b}^Z$ is in $B$ with probability $O(\frac{1}{q}\deg B)$.
Since $\deg B$ is independent of~$q$, this probability is $O(1/q)$.

Similarly, since $\V(f_1,\dotsc,f_Z)$ is of codimension $Z$ for generic polynomials $f_1,\dotsc,f_Z$,
it follows that $\V(f_1,\dotsc,f_Z)$ is of smaller codimension with probability $O(1/q)$ for random polynomials $f_1,\dotsc,f_Z$.
Furthermore, since $\V(f_1,\dotsc,f_Z)$ is defined by $Z$ polynomials, no component
of it can have codimension more than $Z$ (see \cite[Corollary 1.14]{shafarevich}), and so $\V(f_1,\dotsc,f_Z)$
is of pure dimension $b-Z$.

Let $Y=\P^b(\Fq)$. 
Applying \Cref{part:expectproj} we see that
\[
  \Pr\bigl[\abs{\V(f_1,\dotsc,f_Z) \cap \P^b(\Fq)}\leq \tfrac{1}{2}q^{b-Z}\bigr]\leq 4q^{Z-b}.
\]
Since $b-Z\geq 1$, this is also $O(1/q)$, and so random polynomials $f_1,\dotsc,f_Z$ satisfy the conclusion of the lemma with
probability~$1-O(1/q)$.
\end{proof}

\section{Upper bound on \texorpdfstring{$\vd_t(b,m)$}{tau\textunderscore t(b,m)}}\label{sec:vdt}
For the purpose of proving an exponential bound in \Cref{thm:turan}, we need only a bound of the form $\vd_t(b,m)\leq (1-\veps)bt$ that is valid for small~$t$.
We go a step further: it can be shown that $m+2$ points are $m$-dependent if and only if they are collinear, and so $\vd_{m+2}(b,m)=2(b-1)+m$.
The bound on $\vd_t(b,m)$ that we prove is sufficiently strong that the maximum of the quantity $\frac{1}{t-1}\vd_t(b,m)$ in our application
is achieved at $t=m+2$, and so further improvements in the bound would not lead to a smaller base of exponent in \Cref{thm:turan}.

\begin{definition}
  We say that points $p_1,\dotsc,p_t\in \P^b$ are \emph{strongly $m$-dependent} if
  they span $\P^b$ and the associated linear forms $\ell_i(x)=\langle x,p_i\rangle$
  satisfy a relation of the form
  \begin{equation}
    c_1 \prod_{j=1}^m \ell_1(x_j)+\dotsb+c_t \prod_{j=1}^m \ell_t(x_j)=0\label{eq:strongrel}
  \end{equation}
  in which all the coefficients $c_1,\dotsc,c_t$ are non-zero.
\end{definition}
Note that \Cref{eq:hilbequiv} implies that every minimal $m$-dependent set is strongly $m$-dependent inside whichever space it spans.
The key to our upper bound on $\vd_t(b,m)$ is a lower bound on the number of points in any strongly $m$-dependent
set. We begin by proving a relatively weak such bound, adapting the argument which is implicit
in the work of Bia{\l}ynicki\nobreakdash-Birula and Schinzel \cite{bialynicki_schinzel}. We
then bootstrap this weak bound to a stronger bound in \Cref{lem:mdephalf}.

\begin{lemma}\label{lem:weakmdephalf}
  For $m\geq 2$, every strongly $m$-dependent set in $\P^b$ has
  at least $2(b+1)$ many points.
\end{lemma}
\begin{proof}
  Let $p_1,\dotsc,p_t$ be strongly $m$-dependent points in~$\P^b$.
  By renumbering the points if necessary, we may assume that the points $p_1,\dotsc,p_{b+1}$ span $\P^b$.
  If $t\leq 2b+1$, then the remaining points $p_{b+2},\dotsc,p_t$ do not span $\P^b$
  because there are only $t-b\leq b$ of them. Hence, there is an $a\in \P^b$ such that
  $\langle a,p_i\rangle = 0$ for $i=b+2,\dotsc,t$. Plugging $a$ in for each of $x_2, x_3,\dotsc,x_m$, in \eqref{eq:strongrel}
  and writing simply $x$ for $x_1$ we obtain a relation of the form
  \[
    c_1' \ell_1(x)+\dotsb+c_{b+1}' \ell_{b+1}(x)=0
  \]
  where $c_i'=c_i\langle a,p_i\rangle^{m-1}$. Since $p_1,\dotsc,p_{b+1}$ span $\P^b$, not all the coefficients
  $c_i'$ vanish. Because the points $p_1,\dotsc,p_{b+1}$ are linearly independent, the linear forms $\ell_1,\dotsc,\ell_{b+1}$
  are linearly independent, and so we reach a contradiction with the previously-made assumption that $t\leq 2b+1$. 
\end{proof}

To prove the stronger bound, we need to establish a couple of auxiliary results.
\begin{lemma}\label{lem:turann}
  Let $G$ be an $n$-vertex graph with at most $n$ edges. Then $G$ contains an independent set of size
  $n/3$.
\end{lemma}
\begin{proof}
  According to Tur\'an's theorem applied to the complement of $G$, the independence number of $G$ 
  is minimized when $G$ is a union of cliques whose sizes differ by at most~$1$. So, assume that $G$ is of this
  form. This implies that the largest clique in $G$ is of size at most $3$, for otherwise
  $G$ would have had more than $n$ edges. Hence, $G$ is a union of at least $n/3$ disjoint cliques,
  and so has an independent set of size $\lceil n/3\rceil$.
\end{proof}
\begin{lemma}\label{lem:twobases}
  Suppose that $V$ is a finite-dimensional vector space (over any field), $B$ and $B'$ are each a basis for $V$,
  and that no vector in $B$ is a multiple of any vector in $B'$, and vice versa. Then there exists a
  subset $C\subset B$ of size $\abs{C}\geq \tfrac{1}{3}\dim V$ such that $B'\cap \vecspan(C)=\emptyset$.
\end{lemma}
\begin{proof}
  Given a vector $v\in B'$, express it in the basis $B$ as $v=\sum_{b\in B} c_{b,v} b$,
  and define the set $S(v)\eqdef \{b\in B: c_{b,v}\neq 0\}$. Equivalently, the set $S(v)$ is the support
  of $v$, when expressed in the basis~$B$. 

  Let $H\eqdef \{S(v) : v\in B'\}$. Since no vector of $B'$ is a multiple of any vector in $B$, it follows that $\abs{S(v)}\geq 2$,
  and so we may view $H$ as a (potentially non-uniform) hypergraph all whose edges have at least two elements. 
  Since $v\in \vecspan(C)$ if and only if $S(v)\subseteq C$, our task is to find an independent in $H$ of size $n/3$.

  Replace each edge in $H$ by an arbitrary two-element subset thereof, obtaining a graph~$G$. An independent set
  in $G$ is also an independent set in $H$. Since $G$ has $\dim V$ vertices and at most $\dim V$ distinct edges, an 
  appeal to \Cref{lem:turann} completes the proof.
\end{proof}

With this in place, we are ready to prove the stronger bound on~$\vd_t(b,m)$.
\begin{lemma}\label{lem:mdephalf}
  For $m\geq 2$, every strongly $m$-dependent set in $\P^b$ has
  at least $\frac{m+4}{3}(b+1)$ many points.
\end{lemma}
\begin{proof} The proof is by induction on $m$. The base case $m=2$ is contained in \Cref{lem:weakmdephalf}. So,
  assume that $m\geq 3$.

  Let $P$ be strongly $m$-dependent. Let $B\subset P$ be any set of $b+1$ points that span $\P^b$. We break
  into two cases.\smallskip

  \textit{Suppose that the set $P\setminus B$ does not span $\P^b$.} In this case we proceed similarly to the proof of
  \Cref{lem:weakmdephalf}. Namely, we choose an $a\in \P^b$
  such that $\langle a,p\rangle=0$ for all $p\in P\setminus B$. Plugging $a$ in for each of $x_2,x_3,\dotsc,x_m$
  in \eqref{eq:strongrel} and writing $x$ for $x_1$, we obtain a non-trivial linear relation between the $b+1$ linear forms $\langle x,p\rangle$
  with $p\in B$, contradicting the fact that $B$ spans $\P^b$.

  \textit{Suppose that $P\setminus B$ spans $\P^b$.} Let $B'$ be any $b+1$ points of $P\setminus B$ that span
  $\P^b$. Thinking of $B'$ as points in a vector space of dimension $b+1$, from \Cref{lem:twobases} it follows that there is $C\subset B$ of size $\abs{C}\geq (b+1)/3$ such that $B'\cap \vecspan(C)=\emptyset$. Since the field
  $\Fqbar$ is infinite, there is an $a\in \vecspan(C)^{\bot}$ such that $\langle a,p\rangle\neq 0$ for all $p\in P\setminus \vecspan(C)$.
  Setting $x_m=a$ in \eqref{eq:strongrel} hence yields a relation 
  \[
    \sum_{p\in P\setminus \vecspan(C)} c_p' \prod_{j=1}^{m-1} \ell_p(x_j)=0,
  \]
  where $c_p'\eqdef c_p\langle a,p\rangle  \neq 0$. Since the set $P\setminus \vecspan(C)$ contains $B'$, this relation shows that
  $P\setminus \vecspan(C)$ is a strongly $(m-1)$-dependent set. Since
  $\abs{P}\geq \abs{P\setminus \vecspan(C)}+\abs{C}$, we are done by induction.
\end{proof}

Define
\begin{align*}
  \VD'_t(b,m)&\eqdef\{(p_1,\dotsc,p_t)\in \VD_t(b,m) : \text{points }p_1,\dotsc,p_t\text{ span }\P^b\},\\
  \vd'_t(b,m)&\eqdef \dim \VD'_t(b,m).
\end{align*}
\begin{lemma}\label{lem:vdtprime}
Suppose that $m\geq 3$. Then $\vd_t'(b,m)\leq (t-b-1)(b+1)$.
\end{lemma}
\begin{proof}
  For notational brevity, let $\mathcal{L}$ denote the vector space of all
  non-zero linear forms in $b+1$ variables. Let $\mathcal{L}_{*}\eqdef \mathcal{L}\setminus\{0\}$ and define
  \begin{align*}
    U&\eqdef \{ (\ell_1,\dotsc,\ell_t)\in \mathcal{L}_{*}^t : \vecspan \{\ell_1,\dotsc,\ell_t\} = \mathcal{L}\},\\
    W_t'&\eqdef\{(\ell_1,\dotsc,\ell_t)\in U : \sum_{i=1}^t c_i\prod_{j=1}^m \ell_i(x_j)=0\text{ for unique }c_1,\dotsc,c_t\neq 0\},\\
    W_t&\eqdef\{(\ell_1,\dotsc,\ell_t)\in U : \sum_{i=1}^t \prod_{j=1}^m\ell_i(x_j)=0\}.
  \end{align*}
  
  Since each point in $\P^b$ corresponds to a $1$-dimensional family of linear forms (differing up to a multiplication by a non-zero scalar),
  a collection of $t$ points corresponds to a $t$-dimensional family of linear forms.
  Using the characterization of $m$-dependence in \Cref{eq:hilbequiv} this implies that $\dim \VD'_t(b,m)+t\leq \dim W_t'$. 
  Since we also have $\dim W_t'\leq \dim W_t+t$, it follows that $\dim \VD'_t(b,m)\leq \dim W_t$. 
  We shall upper bound $\dim W_t$ by bounding the dimension of the tangent space at every point of~$W_t$.

  Let $(\ell_1,\dotsc,\ell_t)\in W_t$. Since $\ell_1,\dotsc,\ell_t$ span $\mathcal{L}$,
  by renumbering the forms, we may assume that the forms $\ell_1,\dotsc,\ell_{b+1}$ span $\mathcal{L}$.
  Furthermore, by applying a linear change of coordinates, we may also assume that $\ell_i(y_0,\dotsc,y_b)=y_{i-1}$ for each $i=1,2,\dotsc,b+1$.
  Then $(\Delta_1,\dotsc,\Delta_t)\in \mathcal{L}^t$ is in the tangent space to $W_t$ at the point
  $(\ell_1,\dotsc,\ell_t)\in W_t$ if and only if
  \begin{equation}\label{eq:tangentspace}
    \sum_{i=1}^{b+1} \Delta_i(x_k)\prod_{j\neq k} x_{j,i-1}+  \sum_{i=b+2}^t \sum_{k=1}^m \Delta_i(x_k)\prod_{j\neq k} \ell_i(x_j)=0.
  \end{equation}
  Think of this condition as a system of linear equations. Each of the unknowns $\Delta_1,\dotsc,\Delta_t\in \mathcal{L}$
  can be thought of as a vector of $b+1$ many scalar unknowns.  There is one linear equation for each monomial in the $x$-variables appearing in \eqref{eq:tangentspace}.
  For $i\leq b+1$, all the monomials appearing as coefficients of $\Delta_i$ are
  a product of $m-1$ variables of the form $x_{j,i-1}$, and a single other variable. Since $m\geq 3$,
  this means the monomials appearing as coefficients of $\Delta_i$, for $i\leq b+1$, are disjoint.
  This means that the system of equations is of rank at least $(b+1)^2$, i.e., 
  the tangent space to $W_t$ at $(\ell_1,\dotsc,\ell_t)$ is of codimension at least~$(b+1)^2$.
  Therefore,
  \[
    \vd'_t(b,m)\leq \dim W_t\leq t(b+1)-(b+1)^2=(t-b-1)(b+1).\qedhere
  \]
\end{proof}

\begin{lemma}\label{lem:mdependentdim}
Suppose that $t,b,m\geq 3$ are integers. 
Then
\begin{enumerate}
\item \label{part:mdep_empty} $\VD_t(b,m)=\emptyset$ if $t\leq m+1$, and
\item \label{part:mdep_gen} $\vd_t(b,m)\leq  \lfloor\tfrac{3}{m+4}t\rfloor(b+1+\tfrac{m-2}{m+4}t) $ if $m+2\leq t\leq b$.
\end{enumerate}
\end{lemma}
\begin{proof}
From Lagrange interpolation it follows that no set of $m+1$ or fewer points is $m$-dependent. Hence $\VD_t(b,m)=\emptyset$ if $t\leq m+1$. 

Let $r'=\lfloor \tfrac{3}{m+4}t-1\rfloor$. By \Cref{lem:mdephalf}, every minimal $m$-dependent set of size $t$ spans a subspace of projective
dimension at most~$r'$.
Let $\Gr(r,\P^b)$ be the Grassmanian of linear subspaces of dimension $r$ in~$\P^b$.
We then have
\begin{align*}
  \vd_t(b,m)&=\dim \VD_t(b,m)\\
  &=\dim \bigcup_{r\leq r'}\{(p_1,\dotsc,p_t,H)\in \VD_t(b,m)\times \Gr(r,\P^b) : p_1,\dotsc,p_t\text{ span }H \},\\ 
  &= \max_{r\leq r'} \bigl(\vd'_t(r,m)+\dim \Gr(r,\P^b)\bigr)\\
  &\leq \max_{r\leq r'} \bigl((t-r-1)(r+1)+(r+1)(b-r)\bigr)\qquad\qquad\text{by \Cref{lem:vdtprime}}\\
  &=\max_{r\leq r'} (t+b-2r-1)(r+1).\\
\intertext{%
Without the restriction on~$r$, the maximum of the quadratic $(t+b-2r-1)(r+1)$ is achieved when $r=\frac{t+b-3}{4}$.
Since $b\geq t$, the value $\frac{t+b-3}{4}$ is larger than $r'$, and so}
  \vd_t(b,m)&\leq (t+b-2r'-1)(r'+1)\leq (t+b+1-\tfrac{6}{m+4}t)\lfloor \tfrac{3}{m+4}t\rfloor.\hskip-6em&&\qedhere
\end{align*}
\end{proof}

\section{Construction for the Tur\'an problem}\label{sec:turan}
\begin{lemma}\label{lem:turan}
  Let $r,s,Z\geq 1$ be integers. Set $b=r+s+Z$. Suppose that $b,m,Z,r$ satisfy the inequalities $\binom{m+1+r}{m}\geq s^2$, $m\geq 3$
  as well as the condition~\eqref{eq:Zcond} for all fields of sufficiently large characteristic. Then, for every sufficiently large $n$ there is
  a graph with $\Theta(n)$ vertices and $\Omega(n^{2-1/s})$ edges that is $K_{s,t}$-free
  for every $t>m^{s+Z}\prod_{i=1}^r M_i(s^2)$.
\end{lemma}
\begin{proof}
  Using Bertrand's postulate, pick a prime $q$ such that $n\leq q^s<2^sn$. Since $n$ is sufficiently large, we may assume that the condition \eqref{eq:Zcond} is satisfied for~$\Fqbar$. 
  
  Let $f_1,\dotsc,f_Z$ be polynomials as in \Cref{lem:fqindep}. Let $\delta_i\eqdef M_{r-i+1}(s^2)$ for $i=1,2,\dotsc,r$.
  Let $h_1,\dotsc,h_r$ and $h_1',\dotsc,h_r'$ be two independent collections of random homogeneous polynomials on~$\P^b$ with $\Fq$\nobreakdash-coefficients
  of degrees $\deg h_i=\deg h_i'=\delta_i$. Let $g$ be a random bihomogeneous $\Fq$-polynomial on $\P^b\times \P^b$ of bidegree~$(m,m)$.
  
  Define
  \begin{align*}
    L_0&\eqdef \V(f_1,\dotsc,f_Z,h_1,\dotsc,h_r)\cap \P^b(\Fq),\\
    R_0&\eqdef \V(f_1,\dotsc,f_Z,h_1',\dotsc,h_r')\cap \P^b(\Fq).
  \end{align*}
  Let $c'=c'(m,s,\deg W)$ be the constant from \Cref{lem:handy} applied with $W=\V(f_1,\dotsc,f_Z)$, and let $C=C(s^2,s,m^{s+Z},\delta_1,\dotsc,\delta_r)$ be the constant from
  \Cref{cut:diff}. Put $c\eqdef \min\bigl(\tfrac{1}{4},c',(5C)^{-1/s}\bigr)$.
  Let $L$ be a set of size
  $\min(cq^s,\abs{L_0})$ chosen canonically from $L_0$ (for example, it could be the initial segment of $L_0$ in some fixed
  ordering of~$\P^b$). Let $R$ be a similarly defined subset of $R_0$ of size $\min(cq^s,\abs{R_0})$.

  Define the bipartite graph $G$ with parts $L$ and $R$ by connecting the pair $(l,r)\in L\times R$ whenever $g(l,r)=0$. 

  \subparagraph{\texorpdfstring{$G$}{G} has \texorpdfstring{$\Theta(q^s)$}{Theta(q**s)} vertices and \texorpdfstring{$\Omega(q^{2s-1})$}{Omega(q**(2s-1))} edges.}

  By \Cref{part:expectproj}, both $L_0$ and $R_0$ have at least 
  $\tfrac{1}{4}q^s$ elements with probability $1-O(1/q^s)$. Since $c\leq \tfrac{1}{4}$, this means that 
  \begin{equation}\label{rightsize}
    \Pr\bigl[\abs{L}=\abs{R}=cq^s\bigr]\geq 1-O(q^{-s}).
  \end{equation}

  Since the polynomial $g$ is independent of $L_0$ and $R_0$, it follows, by \Cref{part:expectedbi} applied with $Y=L\times R$ to the single random polynomial~$g$,
  that the edge set $E(G)=(L\times R) \cap \V(g)$ is of size at least $\abs{L}\abs{R}/2q$ with probability at least
  $1-O(q^{-2s+1})$. In view of \eqref{rightsize}, it then follows that
  \begin{equation}\label{prob:large}
    \Pr[E(G)\geq \tfrac{1}{2}c^2q^{2s-1}]\geq 1-O(q^{-s}).
  \end{equation}

  \subparagraph{\texorpdfstring{$G$}{G} is \texorpdfstring{$K_{s,t}$}{K\textunderscore s,t}-free.}
  Because of the symmetry between the two parts of $G$, it suffices to show that $G$ is very unlikely to contain $K_{s,t}$ with the part of size $s$ embedded into~$L$.
  Since $L\subset \V(f_1,\dotsc,f_Z)$ and $\V(f_1,\dotsc,f_Z)$ is $m$\nobreakdash-independent, the set $L$ is $m$\nobreakdash-independent.
  Let $W\eqdef \V(f_1,\dotsc,f_Z)$. Since $\binom{m+1+r}{m}\geq s^2$,
  \Cref{lem:handy} applies, telling us that with probability $\tfrac{4}{5}$ every variety of the form
  \[
    W_{l_1,\dotsc,l_s}\eqdef \{w\in W : g(l_1,w)=\dotsb=g(l_s,w)=0\},
  \]
  for distinct $l_1,\dotsc,l_s$,
  is of dimension $\dim W-s=r$. Let $\mathcal{W}$ be the set of all varieties of the form $W_{l_1,\dotsc,l_s}$ for distinct~$l_1,\dotsc,l_s\in L$.
  The set $\mathcal{W}$ is random: it depends on the random choice of polynomials $h_1,\dotsc,h_r$ (because $L$ depends on these polynomials),
  and it depends on the polynomial~$g$. Crucially, $\mathcal{W}$ does not depend on the polynomials $h_1',\dotsc,h_r'$.

  So, we may apply \Cref{cut:diff} to each variety in $\mathcal{W}$ and polynomials $h_1',\dotsc,h_r'$.
  Note that $\deg(W_{l_1,\dotsc,l_s})\leq \deg(W)m^s\leq m^{s+Z}$ by B\'ezout's inequality.
  Therefore, combined with the union bound, the \namecref{prop:randomcut} tells us that
  \begin{align*}
    \Pr[\exists &\text{ distinct }l_1,\dotsc,l_s\in L\text{ s.t.\ \!} \dim \bigl(W_{l_1,\dotsc,l_s}\cap \V(h_1',\dotsc,h_r')\bigr)>0]
    \\&\leq \Pr[ \dim(W_{l_1,\dotsc,l_s})\neq r\text{ for some }l_1,\dotsc,l_s]\\&\qquad+
    \Pr\bigl[\abs{L}\neq cq^s\bigr]+(cq^s)^s \cdot C(s^2,s,m^{s+Z}) q^{-s^2}\\&\leq \tfrac{1}{5}+O(q^{-s})+c^s C(s^2,s,m^{s+Z}).
  \end{align*}
  Because the constant $c$ satisfies $c\leq (5C)^{-1/s}$, this probability is at most $\tfrac{2}{5}+O(q^{-s})$. Note that the variety $W_{l_1,\dotsc,l_s}\cap\nobreak \V(h_1',\dotsc,h_r')$
  contains all common neighbors of the vertices $l_1,\dotsc,l_s$ in the graph $G$. By B\'ezout's inequality, the degree of this variety is at most
  $\deg(W)m^{s}\prod_{i=1}^r \delta_i< t$, and hence
  \[
    \Pr[\text{some }s\text{ vertices in }L\text{ have }t\text{ common neighbors in }G]\leq \tfrac{2}{5}+O(q^{-s}).
  \]
  By symmetry, we may derive the same bound with the roles of $L$ and $R$ reversed, and so
  \begin{equation}\label{prob:kstfree}
    \Pr[G\text{ contains }K_{s,t}]\leq \tfrac{4}{5}+O(q^{-s}).
  \end{equation}

  Putting \eqref{rightsize}, \eqref{prob:large} and \eqref{prob:kstfree} together, it follows that graph $G$ has $\Theta(q^s)=\Theta(n)$ vertices, at least $\Omega(q^{2s-1})=\Omega(n^{2-1/s})$ edges,
  and contains no $K_{s,t}$ with probability at least $\tfrac{1}{5}-O(q^{-s})>0$. In particular, such a graph $G$ exists.
\end{proof}
\begin{lemma}\label{lem:factorial}
For every $r\geq 1$ and every $T$, we have $\prod_{k=1}^r M_k(T)\leq T^{1+\log r}r!$.
\end{lemma}
\begin{proof}
Since $\binom{m+k}{m}\geq \bigl(\frac{m+1}{k}\bigr)^k$, the function $M_k$ satisfies
$
  M_k(T)\leq \lfloor kT^{1/k}\rfloor.
$
The lemma then follows from the inequality $1+\tfrac{1}{2}+\dotsb+\tfrac{1}{r}\leq 1+\log r$.  
\end{proof}
\begin{proof}[Proof of \Cref{thm:turan}]
  We may assume that $s\geq 100$, for otherwise the theorem follows from the result of Alon, R\'onyai and Szab\'o \cite{alon_ronyai_szabo} that we mentioned in the introduction.
  Let $m\eqdef 3$, $r\eqdef \lfloor (6s^2)^{1/3}\rfloor$, and $Z\eqdef s+r+3$. 

  Observe that
  \[
    \tfrac{1}{t-1}\lfloor \tfrac{3}{7}t\rfloor\leq
  \begin{cases}
    \tfrac{1}{2}&\text{ for }t=5,6,7,\\
    \tfrac{7}{15}&\text{ for }t=8,9,10,11,\dotsc.\\
  \end{cases}
  \]
  Hence, for $5\leq t\leq 7$, \Cref{part:mdep_gen} implies that
  \begin{align*}
    \tfrac{1}{t-1}\vd_t(b,m)&\leq \tfrac{1}{2}(s+r+Z+1+\tfrac{1}{7}t)\leq \tfrac{1}{2}(s+r+2)+\tfrac{1}{2}Z<Z.\\
  \intertext{Similarly, for $8\leq t\leq s$ we have}
  \tfrac{1}{t-1}\vd_t(b,m)&\leq \tfrac{7}{15}(s+r+Z+1+\tfrac{1}{7}t)\leq \tfrac{7}{15}(\tfrac{8}{7}s+r+1)+\tfrac{7}{15}Z\\
                          &<\tfrac{8}{15}(s+r+2)+\tfrac{7}{15}Z<Z.
  \end{align*}
  This shows that the condition \eqref{eq:Zcond} holds for $t=5,6,\dotsc,s$.
  It also holds for $t=2,3,4$ by \Cref{part:mdep_empty}.

  We have $\binom{m+1+r}{m}=\binom{r+4}{3}\geq (r+1)^3/6\geq s^2$,
  and so the result follows from \Cref{lem:turan,lem:factorial},
  the estimation
  \[
    r!(s^2)^{1+\log r}\leq 2s^{1/2}(r/e)^r\cdot(s^2)^{1+\log r} \leq s^{\tfrac{2}{3}r+2+2\log r}\leq s^r\qquad\text{for }s\geq 23,
  \]
  and the inequality $3^{r+3}s^r\leq s^{3r/2}\leq s^{4s^{2/3}}$ that is valid for $s\geq 100$.  
\end{proof}

\section{Construction for the Zarankiewicz problem}\label{sec:zar}
This is similar, but simpler than the construction for the Tur\'an problem from the preceding section
because we do not need \Cref{lem:fqindep,lem:mdependentdim}.

\begin{lemma}\label{lem:zar}
  Let $r,s,T\geq 1$ be integers satisfying $T\leq \binom{r+1+m}{m}$. Then, for every sufficiently large $n$ there is
  a sided graph with $\Theta(n^{T/s^2})$ vertices on the left, $\Theta(n)$ vertices on the right, and $\Omega(n^{T/s^2+1-1/s})$ edges that contains no sided $K_{s,t}$
  for every~$t>m^{s}\prod_{i=1}^r M_i(T)$.
\end{lemma}
\begin{proof}
  Let $q$ be the power of $2$ satisfying $n\leq q^s<2^s n$. Let $b\eqdef r+s$.

  Set $\delta_i\eqdef M_{r-i+1}(T)$ for $i=1,\dotsc,r$. Let $c'=c'(m,s,1)$ be the constant from \Cref{lem:handy},
  and let $C=C\bigl(T,s,m^s,\delta_1,\dotsc,\delta_r\bigr)$ be the constant from \Cref{cut:diff}. Put $c\eqdef \min(c',(5C)^{-1/s})$.
  Pick $a\geq 1$ and an $s$-wise $m$-independent set $L$ in $\P^a(\Fq)$ of size $cq^{T/s}$ arbitrarily. For example,
  we may choose $a=cq^{T/s}$ and let $L$ consist of the standard basis vectors.

  We next pick several random polynomials with coefficients in~$\Fq$:
  Let $g$ be a random bihomogeneous polynomial on~$\P^a\times \P^b$ of bidegree~$(m,m)$, and let $h_1',\dotsc,h_r'$ be random
  independently-chosen homogeneous polynomials on~$\P^b$ of degrees $\deg h_i'=\delta_i$. Let
  \[
    R\eqdef \V(h_1',\dotsc,h_r')\cap \P^b(\Fq).
  \]
  Define the sided graph $G$ with the left part $L$ and the right part $R$ by connecting $(l,r)\in L\times R$ whenever $g(l,r)=0$.

  Invoking \Cref{lem:handy} with $W=\P^b$ we see that, with probability $\tfrac{4}{5}$, every variety of the form
  \[
    W_{l_1,\dotsc,l_s}\eqdef \{w\in \P^b: g(l_1,w)=\dotsb=g(l_s,w)=0\},
  \]
  for distinct $l_1,\dotsc,l_s\in L$, has dimension~$b-s=r$. By B\'ezout's inequality, $\deg W_{l_1,\dotsc,l_s}\leq m^s$, and so
  the union bound and \Cref{cut:diff} together imply that
  \begin{equation}\label{zar:kst}
    \Pr[G\text{ contains a sided }K_{s,t}]\leq \tfrac{2}{5}+O(q^{-s}),
  \end{equation}
  where we used that $c\leq (5C)^{-1/s}$, similarly to the corresponding step in the proof of~\Cref{lem:turan}.

  Also, from \Cref{cut:diff} applied to $W=\P^b$ we obtain the inequality $\Pr[\dim \V(h_1',\dotsc,h_r')>s]=\nobreak O(q^{-s})$.
  In view of \Cref{zippel}, this implies that
  \begin{equation}\label{zar:right}
    \Pr\bigl[\abs{R}=O_s(q^s)\bigr]\geq 1-O(q^{-s}).
  \end{equation}
  Finally, by applying \Cref{part:expectproj,part:expectedbi}, it follows that
  \begin{align}
  \notag \Pr\bigl[\abs{E(G)}\!\geq\! \tfrac{1}{4}c q^{T/s+s-1}\bigr]&\mkern-2.1mu\geq\mkern-2.1mu \Pr\bigr[\abs{R}\!\geq\! \tfrac{1}{2}q^s\bigr]\!\Pr\bigl[\abs{E(G)}\!\geq\! \tfrac{1}{2}cq^{T/s}\!\cdot\! \tfrac{1}{2}q^{s-1} \!\bigm|\! \abs{R}\!\geq\! \tfrac{1}{2}q^s\bigr]
  \\\notag &\mkern-2.1mu\geq\mkern-2.1mu \bigl(1-O(q^{-s})\bigr)\bigl(1-O(q^{-s+1}q^{-T/s})\bigr)
  \\&\mkern-2.1mu=\mkern-2.1mu 1-O(q^{-s+1}).\label{zar:edges}
  \end{align}
  
  From \eqref{zar:kst}, \eqref{zar:right}, \eqref{zar:edges} we see that there exists a sided graph
  with $cq^{T/s}$ vertices on the left, $O_s(q^s)$ vertices on the right, and at least $\tfrac{1}{4}cq^{T/s+s-1}$ edges.
\end{proof}
\begin{proof}[Proof of \Cref{thm:zarsmall}]
Let $T\eqdef s^2 \log_n m$, and $r\eqdef \lfloor s/\log^2 s\rfloor $. Note that these constants satisfy $\binom{r+1+k}{r}\geq (s/\log^2 s)^k/k!\geq s^2\log_n m$.
We then use \Cref{lem:zar} with $k$ in place of $m$, and appeal to \Cref{lem:factorial} to bound
\begin{align*}
  \prod_{i=1}^r M_i(T)&\leq T^{1+\log r}r!\leq s^{(1+\log r)k}r^r\\&\leq s^{(1+\log s)\frac{s}{2\log^3 s}}(s/\log^2s)^{s/\log^2 s}\leq e^{2s/\log s}
\end{align*}
to obtain the stated result. 
\end{proof}
\begin{proof}[Proof of \Cref{thm:zarlarge}]
Let $k\geq 1$ be an integer to be chosen shortly.
Let $r\eqdef \lceil s/\log s\rceil$, and $T\eqdef \binom{r+1+k}{k}$
Note that there are constants $c_s''$ and $c_s'''$ such that
\[
  c_s'' k^r\leq T/s^2\leq c_s''' k^r,
\]
and choose $k$ to be the smallest integer so that $c_s''k^r\geq \log_n m$.

Applying \Cref{lem:zar} with $k$ in place of $m$, we obtain a
sided graph 
whose left part is of size $\Omega(n^{T/s^2})=\nobreak \Omega(n^{c_s'' k^r})=\Omega(m)$ and the right part is of size $\Theta(n)$,
matching the K\"ov\'ari, S\'os, Tur\'an bound for $K_{s,t}$-free graphs for
$t>k^s\prod_{i=1}^r M_i(r)$. Since
\[
  k^s\prod_{i=1}^r M_i(r)\leq k^s r^r T^{1+\log r}=O_s(k^{r+2r\log s})=O_s(\log_n m)^{1+2\log s},
\]
this completes the proof.
\end{proof}

\bibliographystyle{plain}
\bibliography{turan}

\end{document}